\documentclass[a4paper,11pt]{article}
\usepackage[utf8]{inputenc}
\usepackage[english]{babel}
\usepackage[T1]{fontenc}
\usepackage[top=3.5cm, bottom=3.5cm, left=3cm, right=3cm]{geometry}
\usepackage{color}    % farver!
\usepackage{graphicx} % pics
\usepackage{listings}
\usepackage{tikz}
\usetikzlibrary{math}
\usepackage{fancyhdr} % layout
\usepackage{lastpage} % lastpage
\usepackage{amsmath,amssymb,mleftright} % amsmath&symboler
\usepackage{listings} % kildekode visning
\usepackage{parcolumns} % minipage
\usepackage{array,tabu,booktabs} % flotte tabeller
\usepackage{multirow} % multirow i tabeller
\usepackage{pgfplots,wrapfig}
\usepackage{titling}  % så vi kan bruge \thedate, \theauthor og \thetitle
\usepackage{lmodern}  % Latin Modern font
\usepackage{verbatim,verbatimbox} % \verbatiminput{<file path>}
\usepackage{fancyvrb} % \begin{Verbatim}[fontsize=\small]
\usepackage{enumitem} % for list: \begin{enumerate}[label=\Alph*]
\usepackage{hyperref} % url
\usepackage{float}
\usepackage{caption}
\usepackage{subcaption}
\usepackage{animate}
\usepackage[final]{pdfpages}
\usepackage{diagbox}
\usepackage{hhline}
\usepackage{titlesec}
\usepackage{algorithm}
\usepackage{algpseudocode}
\usepackage{bigints}
\usepackage[makeroom]{cancel}
\usepackage{xparse}
\usepackage{eufrak}
\usepackage{setspace}
\usepackage{mathtools}
\titlespacing*{\section}
{0pt}{5.5ex plus 1ex minus .2ex}{4.3ex plus .2ex}
\titlespacing*{\subsection}
{0pt}{5.5ex plus 1ex minus .2ex}{4.3ex plus .2ex}

\usepackage[font=small,skip=0pt]{caption}

\usepackage{todonotes} %added by Cilie 30/12

\usepackage{titlesec}
\usepackage{rotating}

\titleformat{\paragraph}
{\normalfont\normalsize\bfseries}{\theparagraph}{1em}{}
\titlespacing*{\paragraph}
{0pt}{3.25ex plus 1ex minus .2ex}{1.5ex plus .2ex}

\makeatletter
\newcommand*\bigcdot{\mathpalette\bigcdot@{.5}}
\newcommand*\bigcdot@[2]{\mathbin{\vcenter{\hbox{\scalebox{#2}{$\m@th#1\bullet$}}}}}
\makeatother

\makeatletter
\def\BState{\State\hskip-\ALG@thistlm}
\makeatother
\definecolor{tempcolor}{rgb}{0.7773 0.0820 0.5195}
\newcommand{\lp}{\left(}
\newcommand{\rp}{\right)}
\newcommand{\lb}{\left[}
\newcommand{\rb}{\right]}

\newcommand{\R}{\mathbb{R}}

\newcommand{\python}{\textsc{Python}}
\newcommand{\fenics}{\textsc{FEniCS}}

\newcommand{\norm}[1]{\ensuremath{\lVert  #1 \rVert}}
\newcommand{\abs}[1]{\ensuremath{\left\vert  #1\right\vert}}

\usepackage{xcolor}
\definecolor{dkgreen}{rgb}{0,0.6,0}
\definecolor{dred}{rgb}{0.545,0,0}
\definecolor{dblue}{rgb}{0,0,0.545}
\definecolor{lgrey}{rgb}{0.9,0.9,0.9}
\definecolor{gray}{rgb}{0.4,0.4,0.4}
\definecolor{darkblue}{rgb}{0.0,0.0,0.6}
\definecolor{turquoise}{rgb}{0.2500,0.8750,0.8125}
\definecolor{indigo}{rgb}{0.2930, 0,0.5078}
\definecolor{mag}{rgb}{1, 0,1}
\definecolor{corn}{rgb}{0.3906,0.5820,0.9258}
\definecolor{mvr}{rgb}{0.7773,0.0820,0.5195}
\definecolor{dod}{rgb}{0.11719,0.5625,1}
\lstdefinelanguage{MatlabCostum}{
      backgroundcolor=\color{white},  
      basicstyle=\footnotesize \ttfamily \color{black} \bfseries,   
      breakatwhitespace=false,       
      breaklines=true,               
      captionpos=b,                   
      commentstyle=\color{dkgreen},   
      emph={repmat, ones},			%Insert missing keywords here
      keywordstyle=\color{blue},           
      escapeinside={\%*}{*)},                  
      frame=single,                  
      language=Matlab,                
      identifierstyle=\color{black},
      stringstyle=\color{blue},      
      numbers=left,                 
      numbersep=5pt,                  
      numberstyle=\tiny\color{black}, 
      rulecolor=\color{black},        
      showspaces=false,               
      showstringspaces=false,        
      showtabs=false,                
      stepnumber=1,                   
      tabsize=4,
      title=\lstname
}
\providecommand{\keywords}[1]
{
  \small	
  \textbf{\textit{Keywords---}} #1
}

\newtheorem{theorem}{Theorem}[section]

\newtheorem{lemma}[theorem]{Lemma}
\newtheorem{proposition}[theorem]{Proposition}
\newtheorem{remark}[theorem]{Remark}

\newcommand{\beginproof}{\emph{Proof.\hspace{4mm}}}%
\newcommand*{\QEDA}{\hfill\ensuremath{\blacksquare}}%
\allowdisplaybreaks

\usepackage[backend=biber,
            style=alphabetic,
            abbreviate=false,
            dateabbrev=false,
            alldates=long]{biblatex}

\begin{document}

%\title{Jacobian of solutions to the conductivity equation in limited view}
%\title{Choice of Neumann boundary values for an elliptic PDE to obtain a non-vanishing Jacobian in limited view }

%\title{Non-vanishing Jacobian for an elliptic PDE in limited view }

\title{Feasibility of hybrid inverse problems in limited view }

\author{Hjørdis Schlüter}
\maketitle

%\tableofcontents
%\newpage

%\section{Abstract}
\begin{abstract}
Hybrid inverse problems such as Acousto-Electric Tomography, Current Density Imaging or Magnetic Resonance Electric Impedance Tomography are concerned with reconstructing the electrical conductivity from interior measurements. For a two-dimensional object the measurements correspond to two different functions imposed as the Neumann boundary condition to an elliptic PDE. In limited view these functions are only non-zero on the part of the boundary that one can control. In this paper we address how to choose such boundary functions in limited view such that the reconstruction procedure is feasible. This is related to the corresponding Jacobian being non-zero. We supplement the theoretical results by numerical reconstructions following the approach of Acousto-Electric Tomography.
\end{abstract}


%\input{Text.tex}
%To dos:
%Fix equation (13)
%Remember to put the real figures in again
%Report the relative error for the reconstruction
%Acknowledgements section
%Fix mathbf notation
%Remark that this also holds for anisotropic conductivities
%Show relative error for both adapted and coordinate functions

\section{Introduction}
In 1980 A.P.Calderón proposed the inverse problem whether it is possible to recover the electrical conductivity by making voltage and current measurements at the boundary. This approach has been realized in medical imaging by the modality \textit{electrical impedance tomography (EIT)}. However, the EIT problem is severely ill-posed and the corresponding reconstructions have low resolution. In order to address these issues, hybrid inverse problems were introduced. For these kind of inverse problems the reconstruction procedure typically consists of two steps: First combine two imaging modalities so that interior measurements can be obtained; secondly recover the conductivity from these interior measurements. We focus on applications that combine the EIT approach with magnetic resonance imaging (MRI) or acoustic waves. The former approaches are known as current density imaging (CDI) \cite{WidlakScherzer, Bal2013_survey, LiSchotlandYangZhong} and Magnetic Resonance Electric Impedance Tomography (MREIT) \cite{Seo_2005,SeoWoo11} and the latter approach is known as Acousto-Electric Tomography (AET) \cite{Zhang:Wang:2004, ammari2008a}. While the EIT approach yields high contrast of the reconstruction, the MRI and acoustic waves yield high resolution of the reconstruction.\\

The mathematical framework behind these applications is the conductivity equation that links the electric potential $u$ with the electrical conductivity $\sigma \in L^{\infty}(\Omega)$:
\begin{equation}\label{condeq}
    \begin{cases}
    -\mathrm{div}(\sigma \nabla u_i)=0 & \text{in }\Omega,\\
    \sigma \nabla u_i \cdot \nu = g_i & \text{on }\partial \Omega,
\end{cases}
\end{equation}
here $\Omega \subseteq \R^2$ models a two-dimensional medium with Lipschitz boundary $\partial \Omega$, $\nu$ is the outward unit normal and the unique solutions $u_i \in H^1(\Omega)$ correspond to different potentials $g_i \in H^{-\frac{1}{2}}(\partial \Omega)$ with $\int_{\partial \Omega} g_i \,\mathrm{d}S=0$ imposed on $\partial \Omega$. 
%The EIT procedure corresponds to recovering $\sigma$ from the potential to current map $$\Lambda_g : g \mapsto \sigma \nabla u \cdot \nu \vert_{\partial \Omega}.$$
The interior measurements for CDI and MREIT are current densities $\sigma \nabla u_i$ corresponding to two different potentials $g_1$ and $g_2$ imposed. For AET the interior measurements are power densities $H_{ij}=\sigma \nabla u_i \cdot \nabla u_j$, with $1 \leq i,j \leq 2$, which are corresponding to two different potentials $g_1$ and $g_2$ imposed. The inverse problem of recovering the conductivity from $\sigma \nabla u_i$ or $\sigma \nabla u_i \cdot \nabla u_j$ requires inversions of the matrix $[\nabla u_1 \, \, \nabla u_2]$. So in order for the reconstruction procedure to be feasible one requires that the following \textit{Jacobian} is non-vanishing:
\begin{equation} \label{eq:nonvan}
    \det [\nabla u_1(x) \, \, \nabla u_2(x)] \neq 0, \quad x \in \Omega.
\end{equation}
The question of how to choose boundary conditions in order to have no vanishing Jacobian has been mainly addressed for Dirichlet boundary conditions. This dates back to Radó, Kneser and Choquet in 1926 and 1945 and is formulated in the Radó-Kneser-Choquet theorem \cite{Rado26,Kneser26,Choquet45}. Here they studied the non-constant coefficient case $\sigma=1$. This result has been generalized to non-constant coefficients and anisotropic conductivities in \cite{Alessandrini86,Alessandrini87,AlessandriniMagnanini94,AlessandriniNesi01,BaumanMariniNesi01,alessandriniNesi15}. For instance, \cite{BaumanMariniNesi01} require that $g=(g_1,g_2)$ is a $C^1$ diffeomorphism and maps $\partial \Omega$ onto the boundary of a convex domain for the condition \eqref{eq:nonvan} to hold. A discussion of such results for Dirichlet boundary conditions is given in \cite{AlbertiCapdeboscq18}. Recent work addresses random boundary functions and the probability for satisfying constraints on the form \eqref{eq:nonvan} \cite{AlbertiArxiv}. \\

We consider a limited view setting, where we only have access to some non-empty closed part of the boundary $\Gamma \subseteq \partial \Omega$ and on the remaining boundary $\sigma \nabla u_i \cdot \nu$ is vanishing:
\begin{equation}\label{condeqlim}
    \begin{cases}
    -\mathrm{div}(\sigma \nabla u_i)=0 & \text{in }\Omega,\\
    \sigma \nabla u_i \cdot \nu = f_i & \text{on }\Gamma,\\
    \sigma \nabla u_i \cdot \nu = 0 & \text{on }\partial \Omega \setminus \Gamma,\\
\end{cases}
\end{equation}
with $f_i \in C(\Gamma)$ and $\int_{\Gamma} f_i \,\mathrm{d}S=0$ so that $\sigma \nabla u_i \cdot \nu \vert_{\partial \Omega} \in H^{-\frac{1}{2}}(\partial \Omega)$ and $u_i \in H^1(\Omega)$. The question of choosing boundary functions $f_1$ and $f_2$ for a non-vanishing Jacobian in limited view has been addressed in \cite{Salo2022} (Dirichlet boundary conditions) and \cite{Jensen2023} (Mixed Dirichlet and Neumann boundary conditions for two- and three-dimensional). We establish similar results to \cite{Salo2022} for this limited 
view setting by extending the result in \cite[Theorem 2.8]{AlessandriniMagnanini94} from no critical points to a non-vanishing Jacobian. Our main result is presented in Theorem \ref{thm:Main} and gives explicit conditions on the boundary functions $f_1$ and $f_2$ in terms of the winding number associated to the curve generated by $f_1$ and $f_2$ (as this not necessarily yields a closed curve we consider a generalized definition). We analyze how boundary functions satisfying these conditions work in practice by reconstructing $\sigma$ numerically from the corresponding power density data $H_{ij}=\sigma \nabla u_i \cdot \nabla u_j$ in limited view following the approach of Acousto-Electric Tomography.  The explicit conditions hold for anisotropic conductivities as well, but we restrict ourself to the isotropic case as this aligns with the numerical examples.\\

Section 2 lists the main result of explicit conditions on the boundary functions $f_1$ and $f_2$. Section 3 lists an analytic reconstruction procedure of $\sigma$ from $\mathbf{H}\in \R^{2\times 2}$ based on \cite{BalMonard}. In section 4 we test boundary functions chosen in accordance with the main result for reconstruction of $\sigma$ from $\mathbf{H}\in \R^{2\times 2}$. We conclude by a discussion in section 5.

\section{Main results}

We state the following proposition without proof corresponding to \cite[Thm 2.8]{AlessandriniMagnanini94}.
\begin{proposition}\label{prop:am94}
    Let $\Omega \subset \R^2$ be a simply connected bounded Lipschitz domain and let $\sigma \in L^{\infty}(\Omega)$ satisfy $\lambda \leq \sigma$ for some $\lambda \in (0,1)$. Let $g \in H^{-\frac{1}{2}}(\partial \Omega)$ with $\int_{\partial \Omega} g \, \mathrm{d}S=0$ and $u \in H^1(\Omega)$ be a nonconstant solution of 
\begin{equation*}
\begin{cases}
    -\mathrm{div}(\sigma \nabla u)=0 & \text{in }\Omega,\\
    \sigma \nabla u \cdot \nu = g & \text{on }\partial \Omega.
\end{cases}
\end{equation*}
Denote by $2M$ the number of closed arcs $\Gamma_1,...,\Gamma_{2M}$ into which $\partial \Omega$ can be split for which $(-1)^j g\geq 0$ on $\Gamma_j$, $j=1,...,2M$. Then the number of interior critical points $N$ satisfies:
\begin{equation*}
    N \leq M - 1.
 \end{equation*}
\end{proposition}

The following lemma is a consequence of Proposition \ref{prop:am94}:

\begin{lemma}\label{lem:4points}
    Let $\Omega, \sigma, u$ and $g$ be as in proposition \ref{prop:am94}. By the proposition it follows that when $\Omega$ has one critical point $x_0 \in \Omega$ so that $\nabla u(x_0)=0$ then $\partial \Omega$ must be composed of at least 4 closed arcs $\Gamma_1, \Gamma_2, \Gamma_3, \Gamma_4$ for which $(-1)^j g \geq 0$. Hence, there must exist points $x_1 \in \Gamma_1, x_2 \in \Gamma_2, x_3 \in \Gamma_3$ and $x_4 \in \Gamma_4$ such that
    \[
\sigma \nabla u \cdot \nu \vert_{x=x_1} < 0, \quad \sigma \nabla u \cdot \nu \vert_{x=x_2} > 0, \quad \sigma \nabla u \cdot \nu \vert_{x=x_3} < 0, \quad \sigma \nabla u \cdot \nu \vert_{x=x_4} > 0.\]
\end{lemma}

We want to revisit the concept of the winding number of a continuous curve $\gamma: [0,\ell] \to \R^2$ and generalize this to the setting were $\gamma$ is not closed, i.e $\gamma(0) \neq \gamma(\ell)$. We assume that $\gamma(t) \neq 0$ for all $t \in [0,\ell]$. We then rewrite $\gamma(t)$ in terms of polar coordinates as 
\[
\gamma(t) = r(t) e^{i\phi(t)}
\]
where $r(t) = \abs{\gamma(t)}$ and $\phi(t)$ are continuous functions in $[0,\ell]$. The angle function 
\[
\mathrm{arg}(\gamma(t)) := \phi(t)
\]
is well defined modulo a constant in $2\pi \mathbb{Z}$. We define 
\[
\mathrm{Ind}(\gamma) := \frac{\mathrm{arg}(\gamma(\ell)) - \mathrm{arg}(\gamma(0))}{2\pi}.
\]
If $\gamma$ is closed, i.e.\ $\gamma(0) = \gamma(\ell)$, then $\mathrm{Ind}(\gamma)$ is the winding number of the curve $\gamma(t)$, which is an integer. If $\gamma$ is not closed, but $\mathrm{arg}(\gamma(t))$ is monotone (i.e.\ nondecreasing or nonincreasing), we may still interpret $\mathrm{Ind}(\gamma)$ as the winding number of $\gamma(t)$, and this is then a real number.

We now give a sufficient condition on a pair of boundary data vanishing outside an arc $\Gamma$ such that the corresponding solutions $u_1, u_2$ satisfy $\mathrm{det} [\nabla u_1(x) \, \nabla u_2(x)]\neq 0$ everywhere in $\Omega$. Note that this can be seen as a consequence of \cite[Theorem 2.8]{AlessandriniMagnanini94}.

\begin{theorem}\label{thm:Main}
Let $\Omega \subset \R^2$ be a bounded simply connected Lipschitz domain with continuous boundary curve $\eta: [0,2\pi] \to \partial \Omega$, and let $\sigma \in L^{\infty}(\Omega)$ satisfy $\lambda \leq \sigma$ for some $\lambda \in (0,1)$. Let $\Gamma = \eta([0,\ell])$ be a closed arc in $\partial \Omega$. Let $f_1, f_2 \in C(\Gamma)$ be linearly independent so that $\int_{\Gamma} f_1 \,\mathrm{d}S=\int_{\Gamma} f_2 \,\mathrm{d}S=0$‚ and assume that $u_i$ is the unique solution of 
\begin{equation}\label{bvpui}
    \begin{cases}
        -\mathrm{div}(\sigma \nabla u_i)=0 & \text{in }\Omega,\\
        \sigma \nabla u_i \cdot \nu=f_i & \text{on }\Gamma ,\\
        \sigma \nabla u_i \cdot \nu = 0 & \text{on }\partial \Omega \backslash \Gamma.
    \end{cases}
\end{equation}
Assume that the curve $\gamma: [0,\ell] \to \R^2$ defined by $\gamma(t) = (f_1(\eta(t)), f_2(\eta(t)))$ satisfies $\gamma(t) \neq 0$ for all $t \in [0,\ell]$, $\mathrm{arg}(\gamma(t))$ is monotone, and $$\abs{\mathrm{Ind}(\gamma)} \leq 1.$$
Then $\mathrm{det} [\nabla u_1(x) \, \nabla u_2(x)]\neq 0$ for all $x \in \Omega$.

\end{theorem}
\beginproof
Assume that the result holds, but one has $\mathrm{det} [\nabla u_1(x_0) \, \nabla u_2(x_0)] =  0$ for some $x_0 \in \Omega$. Then there is a vector $\vec{\alpha} = (\alpha_1, \alpha_2) \in \R^2 \setminus \{0\}$ such that the function 
\[
u = \alpha_1 u_1 + \alpha_2 u_2
\]
satisfies $\nabla u(x_0) = 0$. Note that if $u$ is a constant, then 
\begin{equation*}
    0=\sigma \nabla u \cdot \nu\vert_{\Gamma} = \alpha_1 \sigma \nabla u_1 \cdot \nu\vert_{\Gamma} + \alpha_2 \sigma \nabla u_2 \cdot \nu\vert_{\Gamma} = \alpha_1 f_1 + \alpha_2 f_2 
\end{equation*}
Hence $f_1$ and $f_2$ would be linearly dependent. Thus, we may assume that $u$ is nonconstant. Note that $\sigma \nabla u_i \cdot \nu|_{\partial \Omega}$ are discontinuous, but piecewise continuous and hence they are in $H^{-\frac{1}{2}}(\partial \Omega)$. Then by Lemma \ref{lem:4points} there exist distinct points $x_1, x_2, x_3, x_4 \in \partial \Omega$ that are in counterclockwise order along $\partial \Omega$ such that 
\[
\sigma \nabla u \cdot \nu \vert_{x=x_1} < 0, \quad \sigma \nabla u \cdot \nu \vert_{x=x_2} > 0, \quad \sigma \nabla u \cdot \nu \vert_{x=x_3} < 0, \quad \sigma \nabla u \cdot \nu \vert_{x=x_4} > 0.\]
Consider the function $h: [0,\ell] \to \R$ given by 
\begin{equation}\label{g_fun}
h(t) := \sigma(\eta(t)) \nabla u(\eta(t)) \cdot \nu(\eta(t))= \vec{\alpha} \cdot \gamma(t).
\end{equation}
Extend $h$ by zero to $[0,2\pi)$. Writing $x_j = \eta(t_j)$ where $t_j \in [0,2\pi)$, we have 
\begin{equation} \label{g_contra}
h(t_1) > 0, \quad h(t_2) < 0, \quad h(t_3) > 0, \quad h(t_4) < 0.
\end{equation}
We may assume that $t_1 < t_2 < t_3 < t_4$ (possibly after a cyclic permutation of the indices and after changing $h$ to $-h$).

We now assume that the result holds, and want to derive a contradiction with \eqref{g_contra}. Since $\gamma(t)\neq 0$, $\mathrm{arg}(\gamma(t))$ is monotone and $\abs{\mathrm{Ind}(\gamma(t))} \leq 1$, it follows that $h(t)$ either has at most two zeros in $[0,\ell]$, or has three zeros two of which are at $t=0$ and $t=\ell$. Note that if the argument is not strictly monotone, we make the interpretation that some of these zeros of $h$ could be intervals. Since we assume that $f_1,f_2\in C(\Gamma)$ it follows that $h \in C(\Gamma)$. Now by \eqref{g_contra} $h$ has a zero in the interval $(t_1,t_2)$, a zero in the interval $(t_2,t_3)$ and a zero in the interval $(t_3,t_4)$. Therefore $h$ has at least three zeros in $(0,\ell)$. This is a contradiction. 
\QEDA

\begin{remark}\label{rem:corrdcut}
    Let $\Omega$ be a convex domain and denote by $(f_1^\ell,f_2^\ell)$ normalized cut offs of the coordinate functions $(f_1,f_2)=(x_1,x_2)$. Let $\ell < 2\pi$ so that $\Gamma = \eta([0,\ell])$ is the support of $f_i^\ell$ and $\partial \Omega = \eta([0,2\pi])$ is the support of $f_i$. These normalized cut offs are defined by:
    \begin{equation*}
        f_i^\ell(\eta(t)) = f_i(\eta(t)) \chi_{[0,\ell]}(t) - \frac{\int_{\partial \Omega} f_i(\eta(t)) \chi_{[0,\ell]}(t)\,\mathrm{d}t}{\ell}.
    \end{equation*}
    For convex domains these are automatically captured by Theorem \ref{thm:Main} as $\int_\Gamma f_i^\ell \, \mathrm{d}t=0$ and\\ $\abs{\mathrm{Ind}(f_1(\eta(t)),f_2(\eta(t)))}=1$ so that $\abs{\mathrm{Ind}(f_1^\ell(\eta(t)),f_2^\ell(\eta(t)))}<1$.
\end{remark}

\begin{remark}
For boundary functions $f_1$ and $f_2$ that satisfy one of the conditions in theorem \ref{thm:Main} it is determined by the order of the functions whether $\mathrm{det}[\nabla u_1(x) \, \nabla u_2(x)]$ is positive or negative for all $x \in \Omega$.
\end{remark}

\begin{remark} \label{rem:vio}
Let $\Omega$ be a $C^3$ domain and $f_1,f_2 \in H^\frac{3}{2}(\Gamma)$, then $u_i \in C^{1,\alpha}$ on $\overline{\Omega}$ away from the endpoints of $\Gamma$. Due to the limited view setting $\mathrm{det}[\nabla u_1(x) \, \nabla u_2(x)]$ is zero for $x\in \partial \Omega \setminus \Gamma$ and therefore
\begin{equation*}
    \inf_{x \in \Omega} \mathrm{det}[\nabla u_1(x) \, \nabla u_2(x)] = 0.
\end{equation*}
\end{remark}

\beginproof
The regularity for $u_i$ follows from the fact that if $\sigma \nabla u_i \cdot \nu \vert_{\partial \Omega}(\partial \Omega) \in H^{\frac{3}{2}}(\partial \Omega)$ then $u_i \in H^3(\Omega)$ \cite[Thm 8.29]{salsa2016partial}. This can be extended to this limited view setting so that $u_i \in H^3(\Omega)$ away from the endpoints of $\Gamma$ in the same lines as \cite[Cor. 8.36]{Gilbarg:Trudinger:2001} follows for Dirichlet boundary conditions from \cite[Thm. 8.33]{Gilbarg:Trudinger:2001}. Then by the general Sobolev inequality (see for instance \cite[Thm 6, sec. 5.6.3]{evans98}) it follows that $u_i \in C^{1,\alpha}$ on $\overline{\Omega}$ away from the endpoints of $\Gamma$.
We show that $\mathrm{det}[\nabla u_1(x) \, \nabla u_2(x)]=0$ for the limited view setting by decomposing $\nabla u_i$ into two parts with contribution from the unit outward normal $\nu$ and the unit tangent vector $\omega= \mathcal{J} \nu$ with $\mathcal{J}=\begin{bmatrix} 0 & -1\\ 1 & 0\end{bmatrix}$:
\begin{equation*}
    \nabla u_i = (\nabla u_i \cdot \nu)\nu + (\nabla u_i \cdot \omega)\omega.
\end{equation*}
As $\sigma \nabla u_i \cdot \nu \vert_{\partial \Omega \setminus \Gamma}=0$ by the boundary value problem \eqref{bvpui} it follows that there is no contribution in $\nu$-direction. Therefore both $\nabla u_1\vert_{\partial \Omega \setminus \Gamma}$ and $\nabla u_2\vert_{\partial \Omega \setminus \Gamma}$ are parallel to the unit tangent vector $\omega$ so that $\mathrm{det}[\nabla u_1 \vert_{\partial \Omega \setminus \Gamma} \, \nabla u_2\vert_{\partial \Omega \setminus \Gamma}]=0$. 
\QEDA

\section{Reconstruction procedure}
\label{sec:recproc}
This section lists the analytic reconstruction procedure of an isotropic conductivity $\sigma$ from interior power density measurements $H_{ij}=\sigma \nabla u_i \cdot \nabla u_j$ corresponding two different boundary functions $f_1$ and $f_2$ imposed. This approach is fully based on \cite{monard2012a}. We only list the main equations and refer the reader to \cite{monard2012a} for more details.\\

The reconstruction procedure is split into two steps: Step one to split the functionals $\mathbf{S}_i=\sqrt{\sigma} \nabla u_i$ apart and the second step to recover $\sigma$ from $\mathbf{S}_i=\sqrt{\sigma} \nabla u_i$ for $i=1,2$. 

\subsection{Reconstruction of $\mathbf{S}_i=\sqrt{\sigma} \nabla u_i$}
Reconstructing $\mathbf{S}_i=\sqrt{\sigma} \nabla u_i$ is realized by orthonormalizing the matrix $\mathbf{S}=[\mathbf{S}_1 \, \, \mathbf{S}_2]$ into an $SO(2)$-valued matrix $\mathbf{R}$: $\mathbf{R} = \mathbf{S} \mathbf{T}^T$. The matrix $\mathbf{R}$ is a rotation matrix that is orthogonal and has determinant one and $\mathbf{T}$ is a matrix only depending on the data. $\mathbf{R}$ can be parameterized by an angle function $\theta$:
\begin{equation*}
    \mathbf{R}(\theta) = \begin{bmatrix} \cos \theta & -\sin \theta\\ \sin \theta & \cos \theta \end{bmatrix}.
\end{equation*}
We introduce the functions $\mathbf{V_{ij}}$ that depend on entries of $\mathbf{T}$ and derivatives of its entries:
\begin{equation*}
    \mathbf{V_{ij}}=\nabla (T_{i1})T^{1j}+\nabla (T_{i2})T^{2j}, \quad i,j=1,2
\end{equation*}
here $T^{ij}$ denotes entries in $\mathbf{T}^{-1}$. The first step in the reconstruction procedure then consists of reconstructing $\theta$ and thus $\mathbf{S}$ from the following equation \cite[Eq. (65)]{monard2012a}:
\begin{equation}\label{eq:theta}
    \nabla \theta = \mathbf{F},
\end{equation}
with
\begin{equation*}
    \mathbf{F}=\frac{1}{2}(\mathbf{V_{12}}-\mathbf{V_{21}}-\mathbf{\mathcal{J}}\nabla \log D),
\end{equation*}
where $\mathbf{\mathcal{J}}=\begin{bmatrix} 0 & -1\\ 1 & 0\end{bmatrix}$, and $D=(H_{11}H_{22}-H_{12}^2)^{\frac{1}{2}}$.  Once $\theta$ is known at at least one point on the boundary one can integrate $\mathbf{F}$ along curves originating from that point to obtain $\theta$ throughout the whole domain. If one assumes that $\theta$ is known along the whole boundary one can apply the divergence operator to \eqref{eq:theta} and solve the following Poisson equation with Dirichlet boundary condition:
\begin{equation}\label{eq:thetaBVP}
    \begin{cases}
        \Delta \theta=\nabla \cdot \mathbf{F} & \text{in }\Omega,\\
        \theta=\theta_{\text{true}} & \text{on }\partial \Omega.
    \end{cases}
\end{equation}

\begin{remark}
For implementation we use Gram-Schmidt orthonormalization to obtain the transfer matrix $\mathbf{T}$:
    \begin{equation*}
    \mathbf{T}=\begin{bmatrix} H_{11}^{-\frac{1}{2}} & 0 \\ -H_{12}H_{11}^{-\frac{1}{2}}D^{-1} & H_{11}^{\frac{1}{2}}D^{-1}\end{bmatrix}.
\end{equation*}
This implies that the angle $\theta$ has the interpretation:
\begin{equation*}
    \theta=\text{arg}(\mathbf{R}_1)=\frac{\nabla u_1}{\abs{\nabla u_1}},
\end{equation*}
where $\mathbf{R}_1$ is the first column of $\mathbf{R}$. This implies that $\theta$ describes the angle between $\nabla u_1$ and the $x_1$-axis.
\end{remark}

\subsection{Reconstruction of $\sigma$}
Reconstruction of $\sigma$ is based on \cite[Eq. (68)]{monard2012a}
\begin{equation}\label{eq:sigmaRec}
    \nabla \log \sigma = \mathbf{G},
\end{equation}
with 
\begin{align*}
    \mathbf{G}&=\cos(2 \theta)\mathbf{K} + \sin(2\theta)\mathbf{K},\\
    \mathbf{K}&=\mathcal{U}(\mathbf{V_{11}}-\mathbf{V_{22}})+\mathcal{J}\mathcal{U}(\mathbf{V_{12}}-\mathbf{V_{21}}) \quad \text{and} \quad \mathcal{U}=\begin{bmatrix} 1 & 0\\ 0 & -1\end{bmatrix}.
\end{align*}
Similarly to reconstruction of $\theta$ one can solve this gradient equation by integrating along curves if $\sigma$ is known at one point at the boundary or solving the following boundary value problem if $\sigma$ is assumed known along the whole boundary:
\begin{equation}\label{eq:sigmaBVP}
    \begin{cases}
        \Delta \log (\sigma)= \nabla \cdot \mathbf{G} &\text{in }\Omega,\\
        \log (\sigma)=\log (\sigma_{\text{true}}) & \text{on }\partial \Omega.
    \end{cases}
\end{equation}

\begin{remark}
    For implementation we assume that $\theta$ and $\sigma$ are known along the whole boundary and reconstruct these functions by solving the boundary value problems \eqref{eq:thetaBVP} and \eqref{eq:sigmaBVP}. This requires knowledge of $\theta$ and $\sigma$ at the boundary. Information about the former requires information about $\sigma$ and $u_1$ at the boundary, since these yield information about the tangential component of $\sigma \nabla u_1$. Combining this information with knowledge of $\sigma \nabla u_1 \cdot \nu$ yields information about $\sigma \nabla u_1$ and thus $\theta$. 
\end{remark}

\section{Numerical examples}
\begin{sloppypar}
The \textsc{Python} code to generate the numerical examples can be found on \textsc{GitLab}: \href{https://lab.compute.dtu.dk/hjsc/feasibility-of-hybrid-inverse-problems-in-limited-view}{https://lab.compute.dtu.dk/hjsc/feasibility-of-hybrid-inverse-problems-in-limited-view}.\par
\end{sloppypar}
In this section we want to illustrate how boundary functions can be chosen in accordance with Theorem \ref{thm:Main} and how such functions work in practice to generate suitable AET measurements. We test the "suitability" of the measurements by numerically reconstructing the conductivity from this data following the analytic reconstruction procedure in section \ref{sec:recproc}. For that purpose, we have implemented the reconstruction procedure in section \ref{sec:recproc} in \python{} and use \fenics{}~\cite{fenics} to solve the PDEs. We use a fine mesh to generate our power density data and a coarser mesh to address the reconstruction problem. We use $N_{\text{data}}=79527$ nodes in the high-resolution case, while for the coarser mesh we consider a resolution of $N_{\text{recon}}=50925$ nodes. For both meshes, we use $\mathbb{P}_1$ elements. We consider our domain $\Omega$ to be the unit disk $B(0,1) \subset \R^2$ so that the boundary curve can be parameterized by $\eta:[0,2\pi] \rightarrow \partial \Omega$ as follows:
\begin{equation*}
    \eta(t)=(\cos(t),\sin(t)).
\end{equation*}
We consider 8 limited view settings for $\Gamma_i=\eta([0,\ell_i])$ ranging from $\ell_1=\frac{\pi}{4}$ to the full boundary with $\ell_8=2\pi$:
\begin{equation}\label{eq:Gam}
    \Gamma_i=\eta([0,\ell_i]) \quad \text{with} \quad \ell_i=i \frac{\pi}{4}, \quad \text{for} \quad i=1,..,8.
\end{equation}
We choose two different types of boundary functions for these limited view settings; different types in the sense that one type of boundary functions has more oscillations than the other. The first 8 sets of boundary functions we refer to as \textit{adapted boundary functions} and denote by $(f_{1,\mathrm{adapt}}^{\ell_i},f_{2,\mathrm{adapt}}^{\ell_i})$ in accordance with Theorem \ref{thm:Main} so that $\abs{\mathrm{Ind}(\gamma(t))}=\abs{\mathrm{Ind}(f_{1,\mathrm{adapt}}^{\ell_i},f_{2,\mathrm{adapt}}^{\ell_i})}=1$. These functions are the coordinate functions $(\cos(t),\sin(t)$ adapted to the interval $[0,\ell_i]$ along the boundary:

\begin{equation} \label{eq:fadapt}
(f_{1,\mathrm{adapt}}^{\ell_i},f_{2,\mathrm{adapt}}^{\ell_i})(\eta(t))=\begin{cases}
        (\cos(8t),\sin(8t)) & t \in \lb 0,\ell_1 \rb\\
        (\cos(4t),\sin(4t)) & t \in \lb 0,\ell_2 \rb\\
        \lp \cos \lp \frac{8}{3}t\rp,\sin \lp\frac{8}{3}t\rp \rp & t \in \lb 0,\ell_3 \rb\\
        \lp \cos \lp 2t\rp,\sin \lp 2t\rp \rp & t \in \lb 0,\ell_4 \rb\\
        \lp \cos \lp \frac{8}{5}t\rp,\sin \lp\frac{8}{5}t\rp \rp & t \in \lb 0,\ell_5 \rb\\
        \lp \cos \lp \frac{4}{3}t\rp,\sin \lp\frac{4}{3}t\rp \rp & t \in \lb 0,\ell_6 \rb\\
        \lp \cos \lp \frac{8}{7}t\rp,\sin \lp\frac{8}{7}t\rp \rp & t \in \lb 0,\ell_7 \rb\\
        (\cos(t),\sin(t)) & t \in \lb 0,\ell_8 \rb
    \end{cases}
\end{equation}

These are illustrated in the first column of Figure \ref{fig:funs}. Additionally we choose 8 pairs of boundary functions that corresponds to cut offs of the coordinate functions as in Remark \ref{rem:corrdcut} and therefore we refer to them as \textit{cut off boundary functions} and denote by $(f_{1,\mathrm{cut}}^{\ell_i},f_{2,\mathrm{cut}}^{\ell_i})$. These are so that $\abs{\mathrm{Ind}(\gamma(t))}=\abs{\mathrm{Ind}(f_{1,\mathrm{coord}}^{\ell_i},f_{2,\mathrm{coord}}^{\ell_i})}=\frac{i}{8}$:
\begin{equation}\label{eq:fcoord}
    (f_{1,\mathrm{cut}}^{\ell_i},f_{2,\mathrm{cut}}^{\ell_i})(\eta(t))=\begin{cases}
        (\cos(t),\sin(t))-\lp\frac{2\sqrt{2}}{\pi},\frac{4-2\sqrt{2}}{\pi} \rp & t \in \lb 0,\ell_1 \rb\\
        (\cos(t),\sin(t))-\lp \frac{2}{\pi},\frac{2}{\pi} \rp & t \in \lb 0,\ell_2 \rb\\
        (\cos(t),\sin(t))-\lp-\frac{2\sqrt{2}}{3\pi},\frac{4+2\sqrt{2}}{3\pi} \rp & t \in \lb 0,\ell_3 \rb\\
        (\cos(t),\sin(t))-\lp 0,\frac{2}{\pi} \rp & t \in \lb 0,\ell_4 \rb\\
        (\cos(t),\sin(t))-\lp-\frac{2\sqrt{2}}{5\pi},\frac{4+2\sqrt{2}}{5\pi} \rp & t \in \lb 0,\ell_5 \rb\\
        (\cos(t),\sin(t))-\lp -\frac{2}{3\pi},\frac{2}{3\pi} \rp & t \in \lb 0,\ell_6 \rb\\
        (\cos(t),\sin(t))-\lp-\frac{2\sqrt{2}}{7\pi},\frac{4-2\sqrt{2}}{7\pi} \rp& t \in \lb 0,\ell_7 \rb\\
        (\cos(t),\sin(t)) & t \in \lb 0,\ell_8\rb
    \end{cases}
\end{equation}
These are illustrated in the second column of Figure \ref{fig:funs}.

%\iffalse
\begin{figure}[ht!]
    \centering
    \begin{minipage}[t]{0.5\textwidth}
        \centering
        \includegraphics[width=\textwidth,trim={0cm 1.1cm 0cm 0cm},clip]{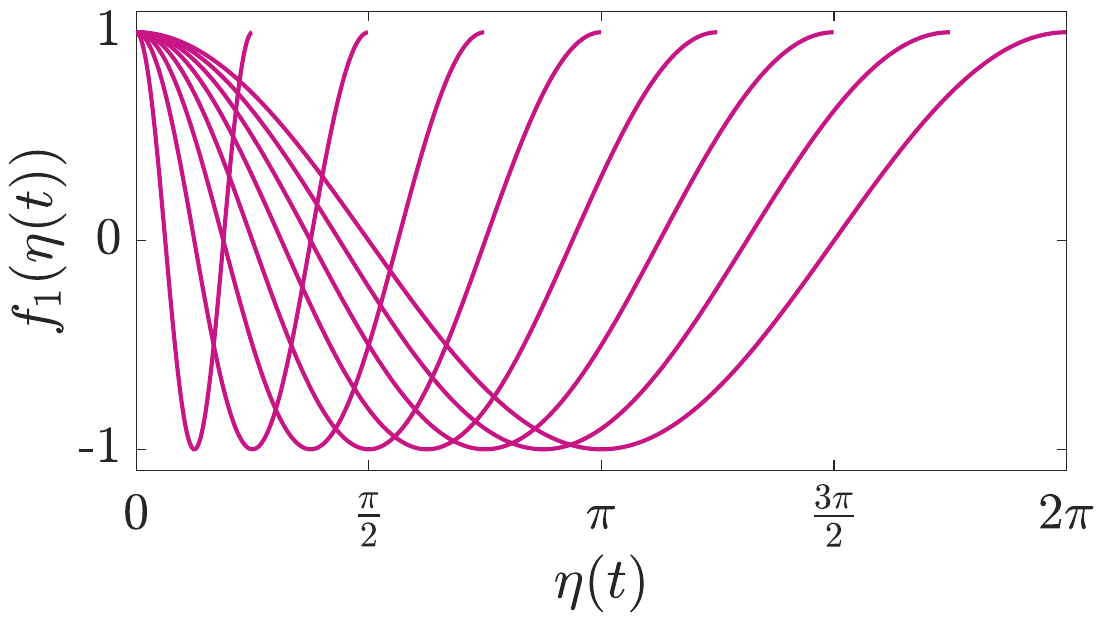}
        \caption*{\hspace{6mm}$t$}
    \end{minipage}%
    \begin{minipage}[t]{0.5\textwidth}
        \centering
        \includegraphics[width=\textwidth,trim={0cm 1.1cm 0cm 0cm},clip]{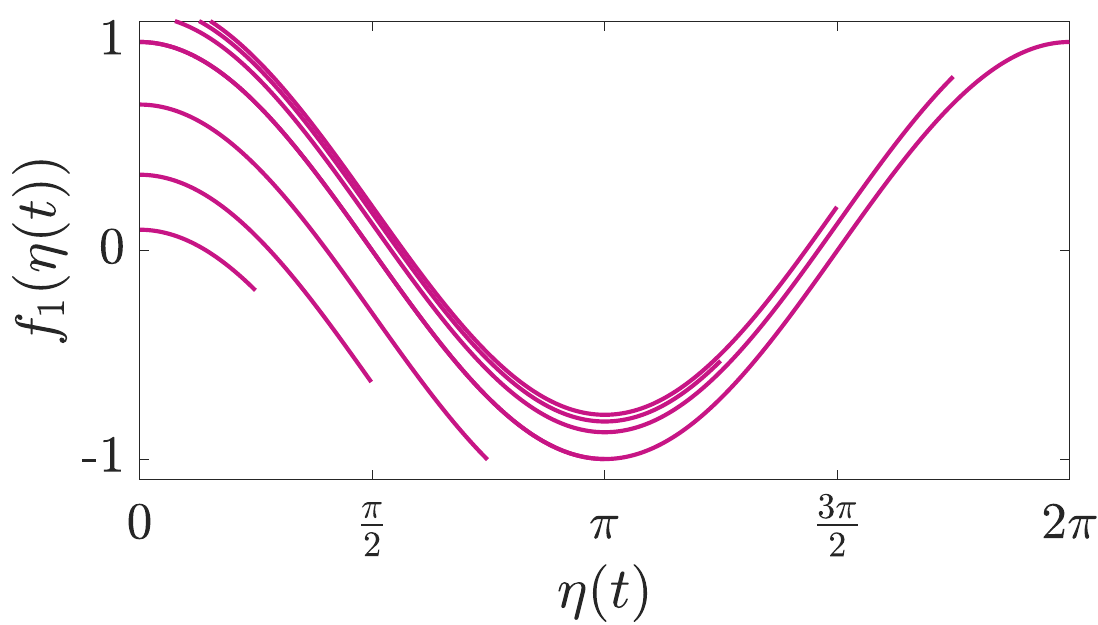}
        \caption*{\hspace{6mm}$t$}
    \end{minipage}
    \begin{minipage}[t]{0.5\textwidth}
        \centering
        \includegraphics[width=\textwidth,trim={0cm 1.1cm 0cm 0cm},clip]{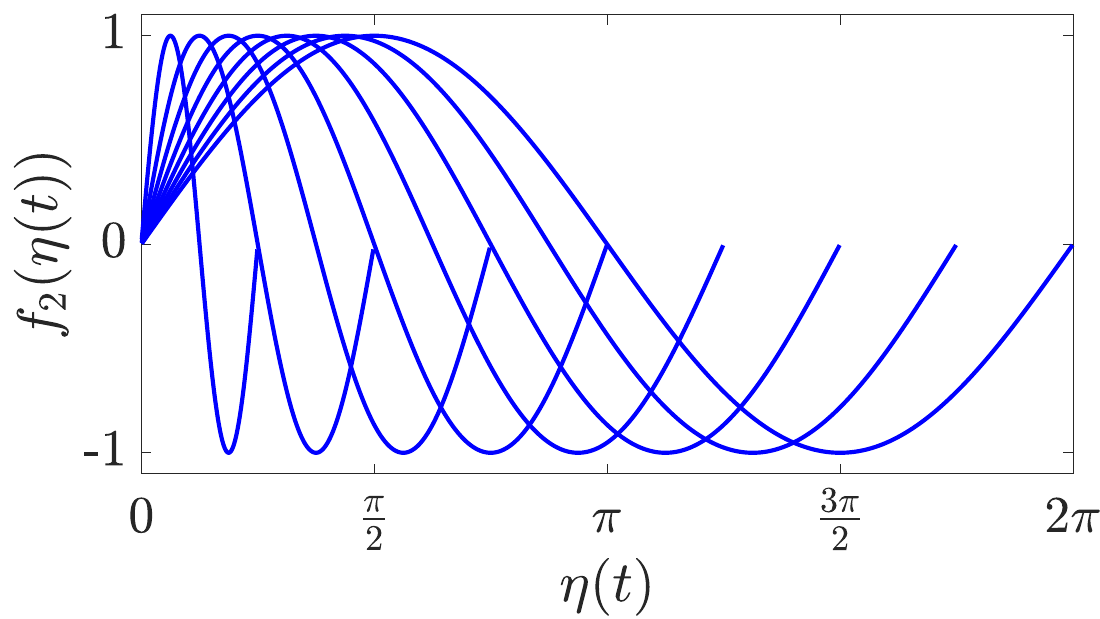}
        \caption*{\hspace{6mm}$t$}
        %\label{fig:f1a}
    \end{minipage}%
    \begin{minipage}[t]{0.5\textwidth}
        \centering
        \includegraphics[width=\textwidth,trim={0cm 1.1cm 0cm 0cm},clip]{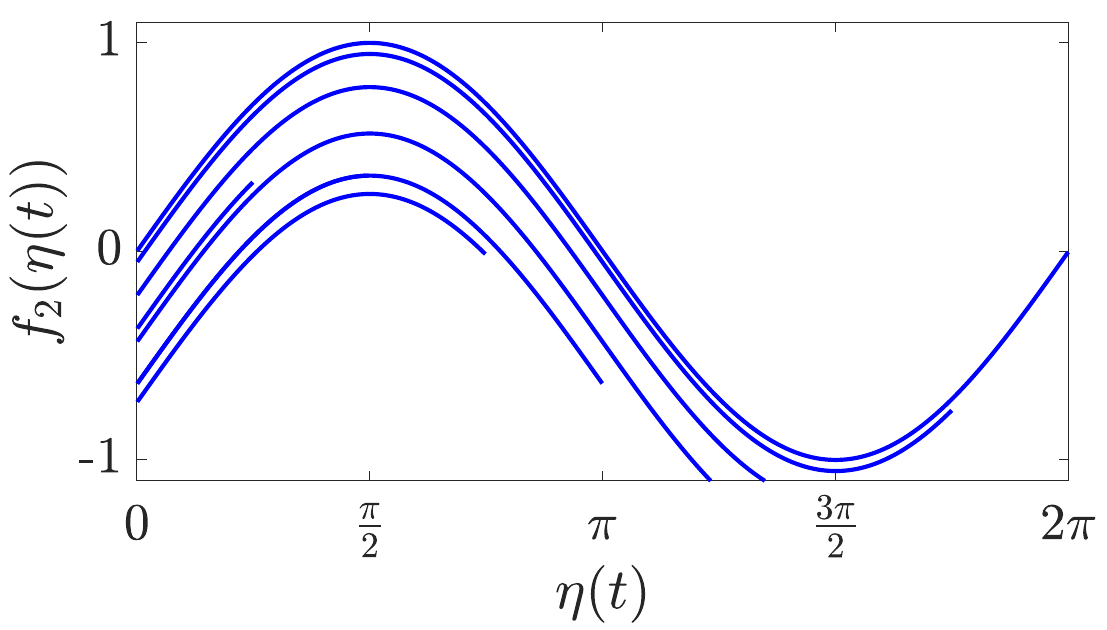}
        \caption*{\hspace{6mm}$t$}
    \end{minipage}
    \caption{Adapted boundary functions as in equation \eqref{eq:fadapt} (left column) and cut off boundary functions as in equation \eqref{eq:fcoord} (right column).}
    \label{fig:funs}
\end{figure}
%\fi

\iffalse
\begin{figure}
    \centering
    \includegraphics[width=\linewidth]{Placeholderfigs/PlaceholdeBfs.png}
    %\caption{Caption}
    \label{fig:funs}
\end{figure}
\fi

We consider two different phantoms; a piecewise constant phantom and a smooth phantom given explicitly by:
\begin{align}\label{eq:sig1}
    \sigma_{\text{1}}(x_1,x_2) &= \begin{cases}2 & \left(x_1+\frac{1}{2} \right)^2 + (x_2)^2 \leq 0.3^2,\\
        2 & \left(x_1\right)^2 + \left(x_2+\frac{1}{2} \right)^2 \leq 0.1^2,\\
        2 & \left(x_1-\frac{1}{2} \right)^2 + \left(x_2-\frac{1}{2} \right)^2 \leq 0.1^2,\\
        1 & \text{otherwise}.
    \end{cases}\\\label{eq:sig2}
    \sigma_{\text{2}}(x_1,x_2) &= \begin{cases}1 + e^{\left(2\,-\,\frac{2}{1-\frac{(x_1)^2+(x_2)^2}{1-0.8^2}}\right)} & 0\leq (x_1)^2+(x_2)^2 \leq 0.8^2,\\
    1 & 0.8^2 \leq (x_1)^2+(x_2)^2 \leq 1,
    \end{cases}
\end{align}
and illustrated in Figure \ref{fig:phantoms}.

%\iffalse
\begin{figure}
    \centering
    \begin{minipage}[t]{0.5\textwidth}
        \centering
        \includegraphics[width=\textwidth]{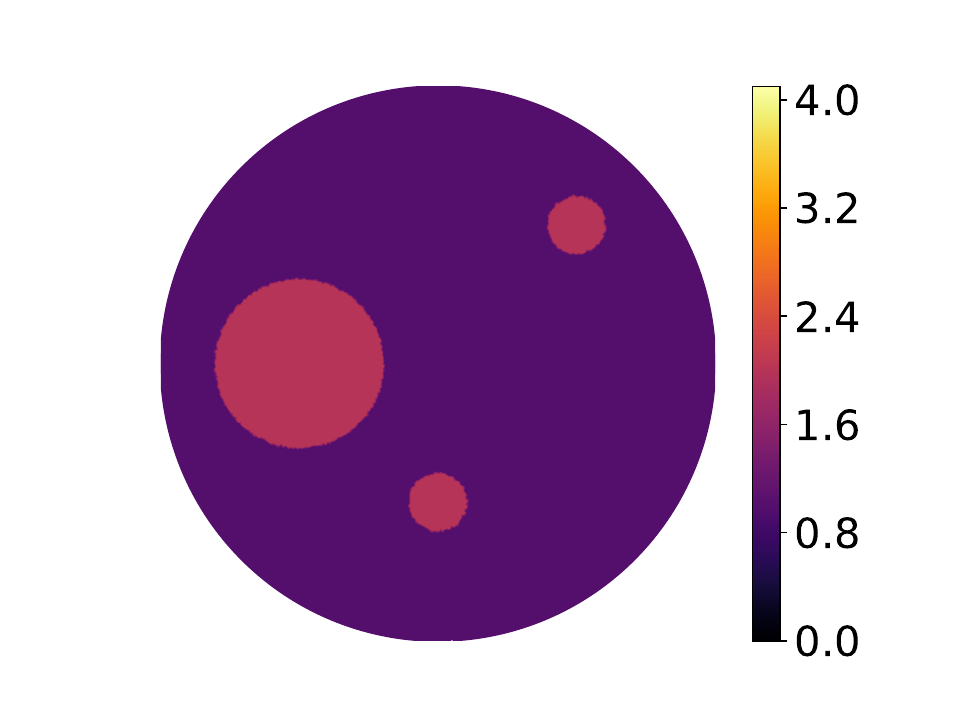}
        \caption*{$\sigma_{\text{1}}$}
    \end{minipage}%
    \begin{minipage}[t]{0.5\textwidth}
        \centering
        \includegraphics[width=\textwidth]{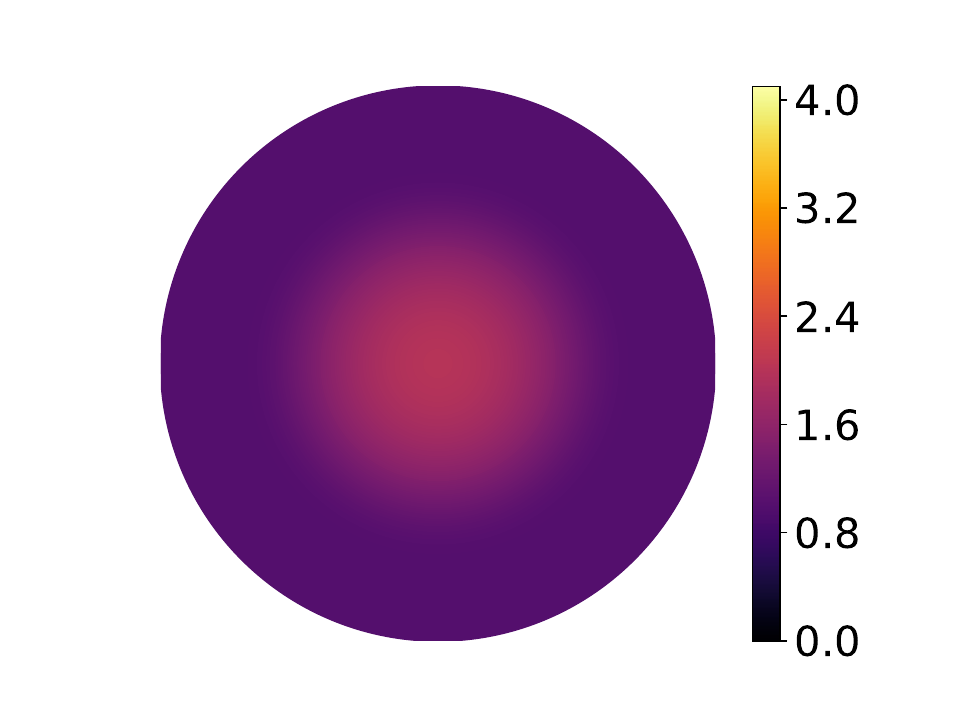}
        \caption*{$\sigma_{\text{2}}$}
    \end{minipage}
    \caption{The true conductivities as in equations \eqref{eq:sig1} and \eqref{eq:sig2}.}
    \label{fig:phantoms}
\end{figure}
%\fi

\iffalse
\begin{figure}
    \centering
    \includegraphics[width=\linewidth]{Placeholderfigs/PlaceholdePhantoms.png}
    %\caption{Caption}
    \label{fig:phantoms}
\end{figure}
\fi

We solve the boundary value problem \eqref{condeqlim} in \fenics{} for the two conductivities in Figure \ref{fig:phantoms}, the 8 limited view settings in equation \eqref{eq:Gam} and the two types of boundary functions in equations \eqref{eq:fadapt} and \eqref{eq:fcoord} to generate power densities. To get an idea of how much information about the conductivity is captured in $H_{11}, H_{12}$ and $H_{22}$ for the two different choices of boundary functions, we illustrate the logarithm of these functions for the limited view setting $\Gamma_4=\eta([0,\pi])$ in Figure \ref{fig:Hijs}. The red line highlights the boundary of control $\Gamma$. As $H_{12}$ holds both positive and negative values we use the following bi-symmetric transformation:
\begin{equation}\label{logsym}
    \log_{\mathrm{sym}}(x) = \mathrm{sign}(x) \log \lp 1+\abs{\frac{x}{C}} \rp
\end{equation}
we use $C=10^{-3}$. This formula is based on \cite[eq.(4)]{webber2013a}. The upper row in Figure \ref{fig:Hijs} shows $H_{11}, H_{12}$ and $H_{22}$ corresponding to the adapted boundary functions in equation \eqref{eq:fadapt} when considering the true conductivity $\sigma_1$ and the limited view setting with $\Gamma_4$. $H_{11}$ and $H_{22}$ visualize how the boundary functions $f_{1,\text{adapt}}$ and $f_{2,\text{adapt}}$ illuminate the domain from different angles and thus capture different bits of information about the conductivity. The cross correlation $H_{12}$ then holds additional information from both boundary functions combined. The edges of the inclusions are clearly visible close to the boundary of control $\Gamma$ and more faded away from $\Gamma$. Similar observations can be made about the second row in Figure \ref{fig:Hijs} corresponding to the cut off boundary functions in equation \eqref{eq:fadapt}. The main difference is that here the illuminated direction in $H_{11}$ for the adapted boundary functions roughly corresponds to the illuminated direction in $H_{22}$ and the same trend for $H_{22}$ for the adapted boundary functions roughly corresponds to the illuminated direction in $H_{11}$. Visually the information stored for both types of boundary functions seems similar. However, as the adapted boundary functions have more oscillations this may play a bigger role for really small boundaries, where these may succeed better at illuminating the domain. 

%\iffalse
\begin{figure}[ht!]
    \centering
    \begin{minipage}[t]{0.33\textwidth}
        \centering
        \input{TikzFigureHijs/AdaptH11}
        %\caption*{1\% noise, \lambda_{\min}=$10^{-7}$}
    \end{minipage}%
    \begin{minipage}[t]{0.33\textwidth}
        \centering
        \input{TikzFigureHijs/AdaptH12}
        %\caption*{5\% noise, \lambda_{\min}=$10^{-6}$}
    \end{minipage}%
    \begin{minipage}[t]{0.33\textwidth}
        \centering
        \input{TikzFigureHijs/AdaptH22}
        %\caption*{10\% noise, \lambda_{\min}=$10^{-6}$}
    \end{minipage}
    \begin{minipage}[t]{0.33\textwidth}
        \centering
        \input{TikzFigureHijs/CoordH11}
        \caption*{$\log(H_{11})$}
    \end{minipage}%
    \begin{minipage}[t]{0.33\textwidth}
        \centering
        \input{TikzFigureHijs/CoordH12}
        \caption*{$\log_{\mathrm{sym}}(H_{12})$}
    \end{minipage}%
    \begin{minipage}[t]{0.33\textwidth}
        \centering
        \input{TikzFigureHijs/CoordH22}
        \caption*{$\log(H_{22})$}
    \end{minipage}
    \caption{Log plots of the power densities in case of adapted boundary functions as in eq. \eqref{eq:fadapt} (upper row) and cut off boundary functions as in eq. \eqref{eq:fcoord} (lower row). The true conductivity is $\sigma_1$ and the limited view setting corresponds to $\Gamma_4$. As $H_{12}$ holds both positive and negative values we use the bi-symmetric transformation $\log_{\mathrm{sym}}$ in \eqref{logsym}.}
    \label{fig:Hijs}
\end{figure}
%\fi

\iffalse
\begin{figure}
    \centering
    \includegraphics[width=\linewidth]{Placeholderfigs/PlaceholdeHijsNew.png}
    %\caption{Caption}
    \label{fig:Hijs}
\end{figure}
\fi

In order to also illustrate the effect of the size of the boundary of control on the information captured in the power densities we show the logarithm of $H_{11}$ corresponding to the adapted boundary functions for the different limited view settings in Figure \ref{fig:H11s}. The edges of the conductivity are clearly visible for the full boundary of control and fade slowly away for decreasing size of $\Gamma$. Until for the smallest boundary $\Gamma_1$ only very little of the small inclusion close to $\Gamma_1$ is visible.

%\iffalse
\begin{figure}[ht!]
    \centering
    \begin{minipage}[t]{0.25\textwidth}
        \centering
        \input{TikzFigureH11s/H11full}
        \caption*{$\Gamma_8$}
    \end{minipage}%
    \begin{minipage}[t]{0.25\textwidth}
        \centering
        \input{TikzFigureH11s/H11Huge}
        \caption*{$\Gamma_7$}
    \end{minipage}%
    \begin{minipage}[t]{0.25\textwidth}
        \centering
        \input{TikzFigureH11s/H11Lar}
        \caption*{$\Gamma_6$}
    \end{minipage}%
    \begin{minipage}[t]{0.25\textwidth}
        \centering
        \input{TikzFigureH11s/H11MedLar}
        \caption*{$\Gamma_5$}
    \end{minipage}
    \begin{minipage}[t]{0.25\textwidth}
        \centering
        \input{TikzFigureH11s/H11Med}
        \caption*{$\Gamma_4$}
    \end{minipage}%
    \begin{minipage}[t]{0.25\textwidth}
        \centering
        \input{TikzFigureH11s/H11SmaMed}
        \caption*{$\Gamma_3$}
    \end{minipage}%
    \begin{minipage}[t]{0.25\textwidth}
        \centering
        \input{TikzFigureH11s/H11Sma}
        \caption*{$\Gamma_2$}
    \end{minipage}%
    \begin{minipage}[t]{0.25\textwidth}
        \centering
        \input{TikzFigureH11s/H11Mini}
        \caption*{$\Gamma_1$}
    \end{minipage}
    \caption{Log plots of the power density $H_{11}$ corresponding to the adapted boundary functions in eq. \eqref{eq:fadapt} and various sizes of the boundary of control $\Gamma$ (indicated by the red line). The true conductivity is $\sigma_1$.}
    \label{fig:H11s}
\end{figure}
%\fi

\iffalse
\begin{figure}
    \centering
    \includegraphics[width=\linewidth]{Placeholderfigs/PlaceholderH11.png}
    %\caption{Caption}
    \label{fig:H11s}
\end{figure}
\fi

In order to check the Jacobian condition \eqref{eq:nonvan} and the influence of the size of the boundary of control we illustrate the logarithm of $\det(\mathbf{H})=\sigma\det([\nabla u_1 \, \nabla u_2])^2$ in Figure \ref{fig:dets}. As the conductivity $\sigma_1$ takes values between 1 and 2, small values of $\det(\mathbf{H})$ correspond to small values of the Jacobian $\det([\nabla u_1 \, \nabla u_2])$ (and thus large negative values of $\log(\det(\mathbf{H}))$ correspond to values of the Jacobian close to zero). For the full boundary of control the small values of the Jacobian are centered around the inclusions and for the limited view settings all smallest values lie towards $\partial \Omega \setminus \Gamma$. This is in accordance with Remark \ref{rem:vio}. For decreasing size of $\Gamma$ the values of $\det(\mathbf{H})$ decrease and especially for the small boundaries of control, $\Gamma_1$ and $\Gamma_2$, almost the whole domain is corrupted by really small values of $\det(\mathbf{H})$ and thus the Jacobian.

%\iffalse
\begin{figure}[ht!]
    \centering
    \begin{minipage}[t]{0.25\textwidth}
        \centering
        \input{TikzFigureDets/DetFull}
        \caption*{$\Gamma_8$}
    \end{minipage}%
    \begin{minipage}[t]{0.25\textwidth}
        \centering
        \input{TikzFigureDets/DetHuge}
        \caption*{$\Gamma_7$}
    \end{minipage}%
    \begin{minipage}[t]{0.25\textwidth}
        \centering
        \input{TikzFigureDets/DetLar}
        \caption*{$\Gamma_6$}
    \end{minipage}%
    \begin{minipage}[t]{0.25\textwidth}
        \centering
        \input{TikzFigureDets/DetLar}
        \caption*{$\Gamma_5$}
    \end{minipage}
    \begin{minipage}[t]{0.25\textwidth}
        \centering
        \input{TikzFigureDets/DetMed}
        \caption*{$\Gamma_4$}
    \end{minipage}%
    \begin{minipage}[t]{0.25\textwidth}
        \centering
        \input{TikzFigureDets/DetSmaMed}
        \caption*{$\Gamma_3$}
    \end{minipage}%
    \begin{minipage}[t]{0.25\textwidth}
        \centering
        \input{TikzFigureDets/DetSma}
        \caption*{$\Gamma_2$}
    \end{minipage}%
    \begin{minipage}[t]{0.25\textwidth}
        \centering
        \input{TikzFigureDets/DetMini}
        \caption*{$\Gamma_1$}
    \end{minipage}
    \caption{Log plots of $\det(\mathbf{H})=\sigma\det([\nabla u_1 \, \nabla u_2])^2$ corresponding to the adapted boundary functions in eq. \eqref{eq:fadapt} and various sizes of the boundary of control $\Gamma$ (indicated by the red line). The true conductivity is $\sigma_1$.}
    \label{fig:dets}
\end{figure}
%\fi

\iffalse
\begin{figure}
    \centering
    \includegraphics[width=\linewidth]{Placeholderfigs/PlaceholderDets.png}
    %\caption{Caption}
    \label{fig:dets}
\end{figure}
\fi

Using the generated power densities we follow the procedure in section \ref{sec:recproc} to reconstruct $\sigma$. %We point out that for some of the limited view settings we observe that the angle $\theta$ changes values from $-\pi$ to $\pi$ throughout $\Omega$; this phenomenon was already observed in \cite{Jensen2023} and \cite{Salo2022} and causes numerical issues as there appears a smoothing of this transition line instead of leaving a discontinuity. For that purpose we add 
We illustrate the reconstructions of $\sigma_1$ using the adapted boundary functions \eqref{eq:fadapt} in Figure \ref{fig:sigmaRecs1} and for $\sigma_2$ in Figure \ref{fig:sigmaRecs2}. The relative errors are shown in Table \ref{tab:sig1} and Table \ref{tab:sig2} respectively. When the boundary of control is large there is visually no big difference in the reconstruction when using either the adapted or the cut off boundary functions. This is also reflected in the relative errors. However, there are some differences in performance for the smaller boundaries of control $\Gamma_1, \Gamma_2$ and $\Gamma_3$. Therefore we compare the reconstructions corresponding to the two types of boundary functions for these smaller boundaries of control in Figure \ref{fig:sigmaRecs2compare} for the conductivity $\sigma_2$. From Figure \ref{fig:sigmaRecs1} and Figure \ref{fig:sigmaRecs1} we see that for the full boundary of control we almost get perfect reconstructions of $\sigma_1$ and $\sigma_2$ as can be seen visually from the reconstructions and from the relative errors of 2.29\% for $\sigma_1$ and 0.03\% for $\sigma_2$. With decreasing size of $\Gamma$ the quality of the reconstruction decreases. Essentially the light areas of the function $H_{11}$ in Figure \ref{fig:H11s} that still capture a lot of information about the conductivity are also the areas for which one gets a good reconstruction. In particular the shape of the dark regions for the small boundaries of control, $\Gamma_1$ and $\Gamma_2$, are exactly the shapes that reappear as shadows in the reconstruction, where the reconstructed values are close to zero and thus way below the actual values between 1 and 2. While the reconstructed values for $\Gamma_1$ and $\Gamma_2$ are too small towards the boundary $\partial \Omega \setminus \Gamma_i$, the reconstructed values for the limited view settings $\Gamma_i$ for $i=3,...,7$ are too large towards this part of the boundary. Apart from the reconstructed values being too small or too large towards $\partial \Omega \setminus \Gamma$ there also appear artifacts close to this part of the boundary. And there seems to happen a transition from $\Gamma_3$ to $\Gamma_2$ where the quality of the reconstruction drops significantly, so this seems to be the point where the boundary of control is too small in order to have enough information for reconstruction when using the procedure in section \ref{sec:recproc}. The reconstruction performance is relatively similar for both conductivities apart from the fact that for the piecewise constant inclusions in $\sigma_1$ the reconstructed value even changes within the inclusions (as can be seen particularly for the boundaries of control $\Gamma_3$ and $\Gamma_4$). For the smooth feature in $\sigma_2$ the reconstructed value is the same for each limited view setting.

%\iffalse
\begin{figure}[ht!]
    \centering
    \begin{minipage}[t]{0.25\textwidth}
        \centering
        \input{TikzFigure/AdaptFull}
        \caption*{$\Gamma_8$}
    \end{minipage}%
    \begin{minipage}[t]{0.25\textwidth}
        \centering
        \input{TikzFigure/AdaptHuge}
        \caption*{$\Gamma_7$}
    \end{minipage}%
    \begin{minipage}[t]{0.25\textwidth}
        \centering
        \input{TikzFigure/AdaptLar}
        \caption*{$\Gamma_6$}
    \end{minipage}%
    \begin{minipage}[t]{0.25\textwidth}
        \centering
        \input{TikzFigure/AdaptMedLar}
        \caption*{$\Gamma_5$}
    \end{minipage}
    \begin{minipage}[t]{0.25\textwidth}
        \centering
        \input{TikzFigure/AdaptMed}
        \caption*{$\Gamma_4$}
    \end{minipage}%
    \begin{minipage}[t]{0.25\textwidth}
        \centering
        \input{TikzFigure/AdaptSmaMed}
        \caption*{$\Gamma_3$}
    \end{minipage}%
    \begin{minipage}[t]{0.25\textwidth}
        \centering
        \input{TikzFigure/AdaptSma}
        \caption*{$\Gamma_2$}
    \end{minipage}%
    \begin{minipage}[t]{0.25\textwidth}
        \centering
        \input{TikzFigure/AdaptMini}
        \caption*{$\Gamma_1$}
    \end{minipage}
    \caption{Reconstructions when using the adapted boundary functions and various sizes of the boundary of control $\Gamma$ (indicated by the red line). The true conductivity is $\sigma_1$ and the colorbar corresponds to the one in Figure \ref{fig:phantoms}.}
    \label{fig:sigmaRecs1}
\end{figure}

\begin{figure}[ht!]
    \centering
    \begin{minipage}[t]{0.25\textwidth}
        \centering
        \input{TikzFigurePhantom2/AdaptFull}
        %\vspace{-5mm}
        \caption*{$\Gamma_8$}
    \end{minipage}%
    \begin{minipage}[t]{0.25\textwidth}
        \centering
        \input{TikzFigurePhantom2/AdaptHuge}
        \caption*{$\Gamma_7$}
    \end{minipage}%
    \begin{minipage}[t]{0.25\textwidth}
        \centering
        \input{TikzFigurePhantom2/AdaptLar}
        \caption*{$\Gamma_6$}
    \end{minipage}%
    \begin{minipage}[t]{0.25\textwidth}
        \centering
        \input{TikzFigurePhantom2/AdaptMedLar}
        \caption*{$\Gamma_5$}
    \end{minipage}
    \begin{minipage}[t]{0.25\textwidth}
        \centering
        \input{TikzFigurePhantom2/AdaptMed}
        \caption*{$\Gamma_4$}
    \end{minipage}%
    \begin{minipage}[t]{0.25\textwidth}
        \centering
        \input{TikzFigurePhantom2/AdaptSmaMed}
        \caption*{$\Gamma_3$}
    \end{minipage}%
    \begin{minipage}[t]{0.25\textwidth}
        \centering
        \input{TikzFigurePhantom2/AdaptSma}
        \caption*{$\Gamma_2$}
    \end{minipage}%
    \begin{minipage}[t]{0.25\textwidth}
        \centering
        \input{TikzFigurePhantom2/AdaptMini}
        \caption*{$\Gamma_1$}
    \end{minipage}
    \caption{Reconstructions when using the adapted boundary functions and various sizes of the boundary of control $\Gamma$ (indicated by the red line). The true conductivity is $\sigma_2$ and the colorbar corresponds to the one in Figure \ref{fig:phantoms}.}
    \label{fig:sigmaRecs2}
\end{figure}
%\fi

\iffalse
\begin{figure}
    \centering
    \includegraphics[width=\linewidth]{Placeholderfigs/PlaceholderRec1.png}
    %\caption{Caption}
    \label{fig:sigmaRec1}
\end{figure}

\begin{figure}
    \centering
    \includegraphics[width=\linewidth]{Placeholderfigs/PlaceholderRec2.png}
    %\caption{Caption}
    \label{fig:sigmaRec2}
\end{figure}
\fi

\begin{table}[ht!]
    \centering
    \caption{Relative $L^2$ errors of $\sigma_1$ for the adapted and cut off boundary functions for varying sizes of $\Gamma$.}
    \begin{tabular}{c||c | c}
         & Adapted functions & Cut off functions\\ \hline \hline
         $\Gamma_8$ & 2.29\% & 2.29\%\\
         $\Gamma_7$ & 21.6\% & 21.2\%\\
         $\Gamma_6$ & 37.5\% & 38.0\% \\
         $\Gamma_5$ & 48.4\% & 49.3\% \\
         $\Gamma_4$ & 52.4\% & 46.1\% \\
         $\Gamma_3$ & 39.9\% & 28.8\% \\
         $\Gamma_2$ & 43.3\% & 55.5\% \\
         $\Gamma_1$ & 84.8\% & 87.1\% \\ \hline
    \end{tabular}
    \label{tab:sig1}
\end{table}

\begin{table}[ht!]
    \centering
    \caption{Relative $L^2$ errors of $\sigma_2$ for the adapted and cut off boundary functions for varying sizes of $\Gamma$.}
    \begin{tabular}{c||c | c}
         & Adapted functions & Cut off functions\\ \hline \hline
         $\Gamma_8$ & 0.03\% & 0.03\%\\
         $\Gamma_7$ & 22.1 \% & 21.2\%\\
         $\Gamma_6$ & 39.1 \% & 38.6\% \\
         $\Gamma_5$ & 51.0\% & 50.5\% \\
         $\Gamma_4$ & 54.0\% & 47.0\% \\
         $\Gamma_3$ & 39.5\% & 26.6\% \\
         $\Gamma_2$ & 37.3\% & 49.9\% \\
         $\Gamma_1$ & 76.2\% & 86.5\% \\ \hline
    \end{tabular}
    \label{tab:sig2}
\end{table}

In Table \ref{tab:sig1} and Table \ref{tab:sig2} we see that for the boundaries of control $\Gamma_1, \Gamma_2$ and $\Gamma_3$ the relative errors differ by more than 10\% for the two choices of boundary functions. This is also reflected visually in Figure \ref{fig:sigmaRecs2compare} when comparing the reconstructions. For $\Gamma_3$ the cut off boundary functions work better, as the value of the reconstructed smooth feature is closer to the original conductivity. For the adapted boundary functions the reconstructed value is too large. For $\Gamma_2$ and $\Gamma_1$ the adapted boundary functions work better, as more of the smooth feature is preserved. For the cut off boundary functions the feature is less visible and the reconstructed values are way too low. A possible explanation for this behavior could be that for $\Gamma_3$ and using the adapted boundary functions the artifacts towards $\partial \Omega \setminus \Gamma_3$ contribute to the reconstructed value of the smooth feature to be too high. For the the smaller boundaries however, all possible information is important so that the more oscillating adapted boundary functions will capture more information, which is lost for the cut off boundary functions.

%\iffalse
\begin{figure}[ht!]
    \centering
    \begin{minipage}[t]{0.33\textwidth}
        \centering
        \input{TikzFigureComparison/Asmamed}  
    \end{minipage}%
    \begin{minipage}[t]{0.33\textwidth}
        \centering
        \input{TikzFigureComparison/Asmall}  
        %\caption*{5\% noise, \lambda_{\min}=$10^{-6}$}
    \end{minipage}%
    \begin{minipage}[t]{0.33\textwidth}
        \centering
        \input{TikzFigureComparison/Amini}      %\includegraphics[width=\linewidth,trim={2.8cm 1cm 2.5cm 1cm},clip]{ComparisonFigures/AdaptMini2.pdf}
        %\caption*{10\% noise, \lambda_{\min}=$10^{-6}$}
    \end{minipage}
    \begin{minipage}[t]{0.33\textwidth}
        \centering
        \input{TikzFigureComparison/Csmamed}  
        \caption*{$\Gamma_3$}
    \end{minipage}%
    \begin{minipage}[t]{0.33\textwidth}
        \centering
        \input{TikzFigureComparison/Csmall}  
        \caption*{$\Gamma_2$}
    \end{minipage}%
    \begin{minipage}[t]{0.33\textwidth}
        \centering
        \input{TikzFigureComparison/Cmini}  
        \caption*{$\Gamma_1$}
    \end{minipage}
    \caption{Comparing reconstructions of $\sigma_2$ when using the adapted boundary functions (upper row) and the cut off of the coordinate functions (lower row) for three different sizes of $\Gamma$ (indicated by red line).}
    \label{fig:sigmaRecs2compare}
\end{figure}
%\fi

\iffalse
\begin{figure}
    \centering
    \includegraphics[width=\linewidth]{Placeholderfigs/PlaceholderComparison.png}
    %\caption{Caption}
    \label{fig:sigmaRecs2compare}
\end{figure}
\fi

\subsection{Reconstruction from noisy measurements}
We want to investigate how noise affects the reconstruction performance for the two choices of boundary functions in limited view. For that purpose we perturb each entry of the power density matrix with normally distributed noise:  
\begin{equation*}
    \widetilde{H}_{ij} = H_{ij} + \frac{\alpha}{100} \frac{e_{ij}}{\norm{e_{ij}}_{L^2}} H_{ij}, \quad 1 \leq i,j \leq 2
\end{equation*}
where $H_{ij}$ is the unperturbed entry, $\widetilde{H}_{ij}$ is the perturbed entry, $\alpha$ is the noise level and $e_{ij}$ is a normally distributed noise so that $e_{ij}\sim \mathcal{N}(0,1)$. We use \texttt{numpy.random.randn} to generate the elements $e_{ij}$ and fix the seed \texttt{numpy.random.seed(50)}. As we assume the power density matrix to be symmetric, we ensure symmetry by computing the perturbed power density matrix by $\frac{1}{2}(\tilde{\mathbf{H}}+\tilde{\mathbf{H}}^T)$. As observed in Figure \ref{fig:dets} the determinant of $\mathbf{H}$ tends to be really small towards the boundary $\partial\Omega \setminus \Gamma$ for the limited view settings relative to the full view setting. Therefore the perturbation above yields negative values of $\det(\tilde{\mathbf{H}})$. In order to ensure that $\mathbf{H}$ is positive definite and thus has a positive determinant we use a lower bound $L$ for the eigenvalues of $\tilde{\mathbf{H}}$.\par
We fix the lower bound $L=10^{-6}$ and consider three different noise levels $\alpha=1, \alpha=5$ and $\alpha=10$. We then compare performance for the two choices of boundary functions and the limited view setting $\Gamma_5=\eta \lp \lb 0, \frac{5\pi}{4} \rb \rp$ in Figure \ref{fig:sigmaRecsNoise5} and the limited view $\Gamma_4=\eta([0,\pi])$ setting in Figure \ref{fig:sigmaRecsNoise4}. The corresponding relative errors are listed in Table \ref{tab:sigNoise5} and Table \ref{tab:sigNoise4}. From the figures we see that the artifacts barely visible in the reconstructions from noise free data are amplified for increasing noise level. The largest artifacts appear close to the end points of $\Gamma$ and the worst artifact appears when using the cut off boundary functions for the limited view setting $\Gamma_5$ at the boundary point $\eta \lp \frac{5\pi}{4}\rp$. This artifact is amplified to such an extent that it dominates the reconstruction for $\alpha=10$. The appearance of the artifacts at the endpoints of $\Gamma$ is explained by the discontinuities in the boundary functions. For the adapted boundary functions $\sigma \nabla u_2 \cdot \nu\vert_{\partial \Omega}$ is continuous, while $\sigma \nabla u_1 \cdot \nu\vert_{\partial \Omega}$ is not. For the cut off boundary functions both $\sigma \nabla u_2 \cdot \nu\vert_{\partial \Omega}$ and $\sigma \nabla u_1 \cdot \nu\vert_{\partial \Omega}$ are discontinuous for the limited view settings and this may induce artifacts close to the discontinuities. However, this does not explain why it is particularly the discontinuity at $t=\frac{5\pi}{4}$ and not at $t=0$ that induces artifacts. We want to point out that only for $\Gamma_5$ there is such a significant difference in performance; for $\Gamma_4$ the performance for both boundary functions is similar. And in this case for both boundary functions artifacts are induced at both endpoints of $\Gamma_4$, but they are not amplified as much for increasing noise level.

\begin{figure}
    \centering
    \includegraphics[width=0.5\textwidth]{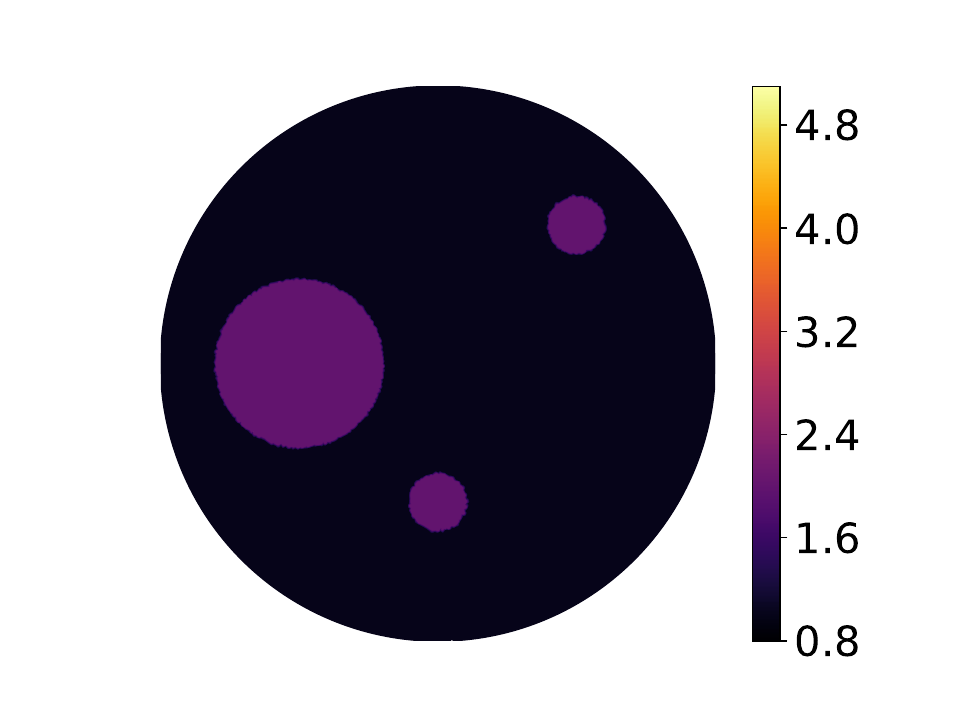}
    \caption{True conductivity $\sigma_1$ used for reconstruction from noisy measurements (note the difference in the range of the colorbar relative to Figure \ref{fig:phantoms}).}
    \label{fig:enter-label}
\end{figure}

%\begin{figure}[ht!]
%    \centering
%        \includegraphics[width=\linewidth]{Figures2/Noisefig.png}
%    \caption{Reconstructions for varying noise levels when using the adapted boundary functions (upper row) and the cut off of the coordinate functions (lower row).}
%    \label{fig:sigmaRecs}
%\end{figure}

\begin{figure}[ht!]
    \centering
    \begin{minipage}[t]{0.33\textwidth}
        \centering
        \input{TikzFigureNoiseMedLar/Adapt1noise}
        %\caption*{1\% noise, \lambda_{\min}=$10^{-7}$}
    \end{minipage}%
    \begin{minipage}[t]{0.33\textwidth}
        \centering
        \input{TikzFigureNoiseMedLar/Adapt5noise}
        %\caption*{5\% noise, \lambda_{\min}=$10^{-6}$}
    \end{minipage}%
    \begin{minipage}[t]{0.33\textwidth}
        \centering
        \input{TikzFigureNoiseMedLar/Adapt10noise}
        %\caption*{10\% noise, \lambda_{\min}=$10^{-6}$}
    \end{minipage}
    \begin{minipage}[t]{0.33\textwidth}
        \centering
        \input{TikzFigureNoiseMedLar/Coord1noise}
        \caption*{$1\%$ noise}
    \end{minipage}%
    \begin{minipage}[t]{0.33\textwidth}
        \centering
        \input{TikzFigureNoiseMedLar/Coord5noise}
        \caption*{$5\%$ noise}
    \end{minipage}%
    \begin{minipage}[t]{0.33\textwidth}
        \centering
        \input{TikzFigureNoiseMedLar/Coord10noise}
        \caption*{$10\%$ noise}
    \end{minipage}
    \caption{Reconstructions of $\sigma_1$ when adding varying noise levels when using the adapted boundary functions (upper row) and the cut off of the coordinate functions (lower row). We consider the limited view setting $\Gamma_5=\eta\lp \lb 0,\frac{5\pi}{4} \rb \rp$ and use the lower bound $L=10^{-6}$ for the eigenvalues of $\tilde{\mathbf{H}}$.}
    \label{fig:sigmaRecsNoise5}
\end{figure}

\begin{figure}[ht!]
    \centering
    \begin{minipage}[t]{0.33\textwidth}
        \centering
        \input{TikzFigureNoise/Adapt1noise}
        %\caption*{1\% noise, \lambda_{\min}=$10^{-7}$}
    \end{minipage}%
    \begin{minipage}[t]{0.33\textwidth}
        \centering
        \input{TikzFigureNoise/Adapt5noise}
        %\caption*{5\% noise, \lambda_{\min}=$10^{-6}$}
    \end{minipage}%
    \begin{minipage}[t]{0.33\textwidth}
        \centering
        \input{TikzFigureNoise/Adapt10noise}
        %\caption*{10\% noise, \lambda_{\min}=$10^{-6}$}
    \end{minipage}
    \begin{minipage}[t]{0.33\textwidth}
        \centering
        \input{TikzFigureNoise/Coord1noise}
        \caption*{$1\%$ noise}
    \end{minipage}%
    \begin{minipage}[t]{0.33\textwidth}
        \centering
        \input{TikzFigureNoise/Coord5noise}
        \caption*{$5\%$ noise}
    \end{minipage}%
    \begin{minipage}[t]{0.33\textwidth}
        \centering
        \input{TikzFigureNoise/Coord10noise}
        \caption*{$10\%$ noise}
    \end{minipage}
    \caption{Reconstructions of $\sigma_1$ when adding varying noise levels when using the adapted boundary functions (upper row) and the cut off of the coordinate functions (lower row). We consider the limited view setting $\Gamma_4=\eta([0,\pi])$ and use the lower bound $L=10^{-6}$ for the eigenvalues of $\tilde{\mathbf{H}}$.}
    \label{fig:sigmaRecsNoise4}
\end{figure}

\iffalse
\begin{figure}
    \centering
    \includegraphics[width=\linewidth]{Placeholderfigs/PlaceholderNoisefig.png}
    %\caption{Caption}
    \label{fig:sigmaRecsNoise}
\end{figure}
\fi

\begin{table}[ht!]
    \centering
    \caption{Relative $L^2$ errors of $\sigma_1$ for the adapted and cut off boundary functions for varying noise levels. The lower bound is $L=10^{-6}$ and the limited view setting is $\Gamma_5=\eta\lp \lb 0,\frac{5\pi}{4} \rb \rp$.}
    \begin{tabular}{c||c | c}
          & Adapted functions & Cut off functions\\ \hline \hline
         $\alpha=1$ & 53.2\% & 59.4\%\\
         $\alpha=5$ & 66.2 \% & 235\%\\
         $\alpha=10$ & 94.6 \% & 714\% \\\hline
    \end{tabular}
    \label{tab:sigNoise5}
\end{table}

\begin{table}[ht!]
    \centering
    \caption{Relative $L^2$ errors of $\sigma_1$ for the adapted and cut off boundary functions for varying noise levels. The lower bound is $L=10^{-6}$ and the limited view setting is $\Gamma_4=\eta([0,\pi])$.}
    \begin{tabular}{c||c | c}
          & Adapted functions & Cut off functions\\ \hline \hline
         $\alpha=1$ & 65.1\% & 67.8\%\\
         $\alpha=5$ & 72.9 \% & 70.7\%\\
         $\alpha=10$ & 80.7 \% & 76.0\% \\\hline
    \end{tabular}
    \label{tab:sigNoise4}
\end{table}

%\begin{figure}
%    \centering
%    \includegraphics[width=\textwidth,trim={2.6cm 1.1cm 4cm 1.1cm},clip]{AdaptedNoise/ContMedNoise11e7.pdf}
%    \caption{Caption}
%    \label{fig:enter-label}
%\end{figure}
\clearpage

\section{Discussion}
This paper presented explicit conditions on the Neumann boundary functions so that the corresponding solutions have a non-vanishing Jacobian in limited view. The non-vanishing Jacobian implies that the reconstruction procedure for an isotropic conductivity $\sigma$ from power density measurements $H_{ij}=\sigma \nabla u_i \cdot \nabla u_j$ for $1 \leq i,j \leq 2$ as in Acousto-Electric Tomography or from current density measurements $\sigma \nabla u_i$ for $i=1,2$ as in Current Density Imaging or Magnetic Resonance Electric Impedance Tomography is feasible.\par
We extended the idea for Dirichlet boundary functions in \cite{Salo2022} to Neumann boundary functions by using tools from complex analysis. We showed in Theorem \ref{thm:Main} that the boundary functions should not be too oscillating in the sense that the absolute value of the winding number associated to the curve generated by the boundary functions $f_1$ and $f_2$ should be at most 1 (as this not necessarily yields a closed curve we consider a generalized definition of the winding number). Additionally, in order to guarantee unique solutions to the boundary value problem $f_1$ and $f_2$ should integrate to zero. This is in contrast to \cite{Salo2022}, where the restrictions on the oscillations is formulated in terms of $f_1'$ and $f_2'$ and where one distinguishes between boundary functions that have discontinuities at one or both end points of $\Gamma$. Due to the restrictions on $f_1$ and $f_2$ (for Neumann boundary functions) instead of restrictions on their derivatives (for Dirichlet boundary functions) it does not make a difference in the proof whether the functions extend continuously to zero at one end of $\Gamma$ or they are discontinuous at both ends. In particular, with the restrictions both on the oscillations of $f_1$ and $f_2$ and the functions integrating to zero it is difficult to construct functions for which $\sigma \nabla u_i \cdot \nu \vert_{\partial \Omega}$ is not discontinuous. In fact, it is not possible to choose two linearly independent functions where both are trigonometric polynomials that extend continuously to zero at both ends of $\Gamma$ (only $f_{2,\mathrm{adapt}}$ as in equation \eqref{eq:fadapt} satisfies this criteria). Therefore for the numerical examples we restricted ourselves to $f_1$ and $f_2$ being trigonometric polynomials, where at least one of the functions yields discontinuities in $\sigma \nabla u_i \cdot \nu \vert_{\partial \Omega}$. \par
While the discontinuities of the Neumann boundary functions do not play a role in the proof of Theorem \ref{thm:Main} they do play a role in the numerical reconstruction of $\sigma$ from noisy power density measurements corresponding to boundary functions chosen in accordance with Theorem \ref{thm:Main}. Depending on the choice of $\Gamma$ there appear artifacts towards one or both end points of $\Gamma$ that are amplified for increasing noise level. For reconstructions from noise free data we see little effect of the discontinuities in the boundary functions, but the oscillations of the functions seem to play a role. For limited view settings, where one only has control over a quarter or less of the boundary it seems that more oscillating boundary functions seem to illuminate the domain better so that more properties of the conductivity can be reconstructed. For limited view settings, where the boundary of control is larger, the oscillations seem to be less significant. \par
In possible future work it could be interesting to quantify how far away from the boundary of control one can trust the reconstruction of the conductivity.

\section{Acknowledgments}
H.S. was supported by the Research Council of Finland (Flagship of Advanced Mathematics for Sensing Imaging and Modelling grant 359208; Centre of Excellence of Inverse Modelling and Imaging grant 353092; and other grants 351656, 358047).

\clearpage

%\newpage
%\input{Conclusions.tex}

\clearpage
\newpage
\printbibliography
%\newpage

Department of Mathematics and Statistics, University of Jyväskylä, Finland\\
\indent \textit{E-mail address: } \texttt{hjordis.a.schluter@jyu.fi}
%\texttt{hjsc@dtu.dk}

%\input{Appendix.tex}

\end{document}